\documentclass[a4paper,10pt]{article}
\usepackage[utf8]{inputenc} 
\usepackage{float}

\usepackage[maxfloats=500]{morefloats}
\begin{document}        
\begin{center} 
{\bf{\Large Calculation and Recurrence Relations of the Determinants of the Hessenberg - Toeplitz Matrices
  }}
\vspace{1cm}

E. Kyriakopoulos\\
Department of Physics\\ 
National Technical University
15780 Zografou, Athens, GREECE\\
E-mail: kyriakop@central.ntua.gr\\
\end{center}
\begin{abstract}
The determinant $D(n)$ of a $n*n$  Hessenberg - Toeplitz matrix
is given as product of certain functions of its elements. A compact expression of $D(n)$ is presented and recurrence relations involving determinants and matrix elements of the Hessenberg -Toeplitz matrices are found.
\end{abstract}

Classification numbers: 15B05, 15A99


.

The determinants of the  Hessenberg-Toeplitz matrices are interested quantities and several papers for these determinants have been published \cite{IN1}, \cite{ME1}, \cite{GO1}.   
In this paper we diagonalize the Hessenberg-Toeplitz matrix
\begin{equation}
T(n)={ 
\left|
\begin{array}{cccccc}
 C(1) & C(2) & \dots & C(n-2) & C(n-1) & C(n)\\

1 & C(1) & C(2) & \dots & C(n-2) & C(n-1) \\

0 &  1 & C(1) & \dots &  C(n-3) & C(n-2)  \\

\vdots & \vdots  &  \vdots & \vdots  &  \vdots \\

0  & 0  & 0  & \dots &  C(1)  &  C(2) \\

0  &  0  &  0  &  \dots &  1  &  C(1)\\
 
 \end{array} 
 \right|
}
\label{1}
\end{equation}
 and calculate its determinant $D(n)$. To do that we multiply the first row of the matrix by -$\frac{1}{C(1)}$ and add to the second row. Then we multiply the second row which comes after the addition by -$\frac{1}{  C(1)-
\frac{C(2)}{C(1)}   }$ and add to the third row.Then we multiply the third row which comes after the addition by
-$\frac{1}{ C(1)    -\frac{C(2)- \frac{C(3)}{(C(1)  }}{C(1)-\frac {C(2)}{C(1)}}  }$   and add to the fourth row.  If we continue to work in this way we find a matrix with diagonal elements $d(i)$ $i= 1, 2, 3, ...n$  where    
\begin {eqnarray}
&& d(1)=C(1) \nonumber\\
&& d(2)=C(1)-\frac{C(2)} {C(1)} \nonumber\\
&&  d(3)=C(1)  -\frac{C(2)-\frac{C(3)}{C(1)}}   { C(1)-\frac{C(2)}{C(1) } }  \nonumber\\
&&  d(4)=C(1)-\frac{C(2)  -\frac{C(3)-\frac{C(4)}{C(1)}}   { C(1)-\frac{C(2)}{C(1) } }  }   {C(1)  -\frac{C(2)-\frac{C(3)}{C(1)}}   { C(1)-\frac{C(2)}{C(1) } }  }  \nonumber\\
&& d(5)=C(1)-\frac{d'}{d(4)}\nonumber\\
&&  .   \nonumber\\
&&  .   \nonumber\\
&&   .
\label{2}
\end{eqnarray}
where $d'$ is given by the relation
\begin{eqnarray}
&&  d'=C(2)-\frac{C(3)  -\frac{C(4)-\frac{C(5)}{C(1)}}   { C(1)-\frac{C(2)}{C(1) } }  }   {C(1)  -\frac{C(2)-\frac{C(3)}{C(1)}}   { C(1)-\frac{C(2)}{C(1) } }  } 
\label{3}   
\end{eqnarray}
and all elements below the diagonal elements equal to zero. If we define $F_1(j)$  and $F_i(j)$ by the relations
\begin{equation}
F_1(j)=C(j), \;\;;j=1, 2,3,  
\label{4}
\end{equation}
\begin{equation}
F_i(j)=C(j)-\frac{F_{i-1}(j+1)}{F_{i-1}(1)} , \;\; \;\;  i=2, 3, 4,  ..., ,j=1, 2, 3, ...
\label{5}
\end{equation}
we can show that for any $i$ we have
\begin {equation}
d(i)=F_{i}(1)
\label{6}
\end{equation}
and 
we get the following Theorem 1:

Theorem 1

The determinant $D(n)$ of the Hessenberg-Toeplitz matrix of Eq (\ref{1}) is given by the relation
\begin{equation}
D(n)= F_1(1) F_2(1) F_3(1)... F_{n}(1)
\label{7}
\end{equation}
where  $F_i(j)$ are given by Eqs (\ref{4})  and  (\ref{5}) .

Since the quantities $F_i(j)$ are functions of $C(j)$ Eq (\ref{7})  gives $D(n)$ as a polynomial of $C(j)$. To find this polynomial we write
\begin{eqnarray}
&&  C(1)C(j)-C(j+1)=C(1,j)  \nonumber\\
&&  C(1,1)C(j)-((1,j+1)=C(1,1,j) \nonumber\\
&&  C(1,1,1)C(j)-C(1,1,j+1)=C(1,1,1,j)   \nonumber\\
&&...
\label{8}
\end{eqnarray}
or in general
\begin{equation}
 C(a_1,a_2,...a_n)C(j)-C(a_1,a_2,..a_{n-1},j+1)=C(a_1,a_2,...a_n,j)
 \label{9}
 \end{equation}
 where $a_j=1, j=1, 2, 3,... $
Then we get using  Eqs (\ref{2}),(\ref{4})-(\ref{6}), (\ref{8}) and (\ref{9})
\begin{eqnarray}
    && F_4(1)=d(4)=C(1)-\frac{C(2)  -\frac{C(1)C(3)-C(4)}{C(1)C(1)-C(2)}}{C(1)  -\frac{C(1)C(2)-C(3)}{C(1)C(1)-C(2)}}=C(1) -\frac{C(2)  -\frac{C(1,3)}{C(1,1)}}{C(1)  -\frac{C(1,2)}{C(1,1)}}=  \nonumber\\
    && C(1) - \frac{C(1,1)C(2)-C(1,3)}  
    {C(1,1)C(1) -C(2)}=C(1)-\frac{C(1,1,2)}{C(1,1,1)}=\frac{C(1,1,1,1)}{C(1,1,1)}
\label{10}
\end{eqnarray}
In a similar way we find that
\begin {eqnarray}
&& F_2(1)=\frac{C(1,1)}{C(1)},\;\;\;
F_3(1)=\frac{C(1,1,1)}{C(1,1)},\;\;\; F_5(1)=\frac{C(1,1,1,1,1)}{C(1,1,1,1)}, \;\;\;     \nonumber\\
&& F_6(1)=\frac{C(1,1,1,1,1,1)}{C(1,1,1,1,1)}
\label{11}
\end{eqnarray}
In general we have
\begin{equation}
F_n(1)=\frac{C(a_1,a_2,a_3,...a_{n-1},a_n)}{C(a_1,a_2,a_3,...a_{n-2},a_{n-1})}
\label{12}
\end{equation}
where $a_j=1$ for $j=1, 2. 3. ...n$. Then from Eqs (\ref{7})  and  (\ref{12}) we find the following Theorem 2

Theorem 2

The determinant $D(n)$  of the Hessenberg-Toeplitz matrix of Eq (\ref{1}) for arbitrary $n$ is given by the relation  \begin{equation}
D(n)=C(a_1,a_2,a_3,...a_n)\;\;with \;\;a_j=1\;\;\;j=1,2,3,...n
\label{13}
\end{equation}
where the quantities  $C(a_1,a_2,a_3,...a_n)$ are given by Eq (\ref{9})

We define $D(1)$ by the relation
\begin{equation}
    D(1)=C(1)
    \label{14}
\end{equation}
and we get from  Eqs (\ref {1}) and (\ref {8}) the relation
\begin {equation} 
D(2)=C(1,1)
\label{15}
\end{equation} 
Then we find from  Eqs (\ref {8}), (\ref {9}) and (\ref {13}) - (\ref {15}) the relation
\begin{eqnarray} 
&& D(n)= C(a_1,a_2,...a_n)= C(a_1,a_2,...a_{n-1})C(1)-C(a_1,a_2,...a_{n-2},2)=  \nonumber\\
&& D(n-1) C(1)-C(a_1,a_2,...a_{n-2},2))=D(n-1)C(1)-C(a_1,a_2,...a_{n-2} )C(2)  \nonumber\\
&& +C(a_1,a_2,...a_{n-3} ,3))=D(n-1)C(1)- \nonumber\\
&& D(n-2)C(2)+C(a_1,a_2,...a_{n-3} ,3)=...
\label{16}
\end{eqnarray}
or
\begin{eqnarray} 
&& D(n) = D(n-1)C(1)-
D(n-2)C(2)+D(n-3)C(3)-... \nonumber\\
&&+(-1)^{n-3}D(2)C(n-2)+ (-1)^{n-2}D(1)C(n-1) +(-1)^{n-1}C(n) \nonumber\\
&& n=2, 3, 4,...
\label{17}
\end{eqnarray}
Also  we find the relation
\begin{eqnarray}
&& D(n)= C(a_1,a_2,...a_n)= (C(a_1,a_2,...a_{n-1})C(1)- \nonumber\\
&& C(a_1,a_2,...a_{n-2},2)=(C(a_1,a_2,...a_{n-2})C(1)-C(a_1,a_2,...a_{n-3},2))C(1) \nonumber\\
&& -C(a_1,a_2,...a_{n-2})C(2)+C(a_1,a_2,...a_{n-3},3)=C(a_1,a_2,...a_{n-2}) \nonumber\\
&&  C(1,1)-C(a_1,a_2,...a_{n-3},2))C(1)+C(a_1,a_2,...a_{n-3},3)=D(n-2)C(1,1) \nonumber\\
&&-(C(a_1,a_2,...a_{n-3})C(2)-C(a_1,a_2,...a_{n-4},3))C(1)+C(a_1,a_2,...a_{n-3}) \nonumber\\
&& C(3)-C(a_1,a_2,...a_{n-3},3)=D(n-2)C(1,1)-D(n-3)C(1,2)+ \nonumber\\
&& C(a_1,a_2,...a_{n-4},3))C(1)-C(a_1,a_2,...a_{n-3},4)=...
\label{18}
\end{eqnarray}
or
\begin{eqnarray} 
&& D(n) = D(n-2)C(1,1)-
D(n-3)C(1,2)+D(n-4)C(1,3)-... \nonumber\\
&&+(-1)^{n-3}D(2)C(1,n-2)+ (-1)^{n-1}D(1)C(1,n-2) + \nonumber\\
&& (-1)^nC(1,n-1)\;\;\;n=3, 4, 5,...
\label{19}
\end{eqnarray}

Proceeding in the same way we find  Theorem 3

Theorem 3

In addition to the recurrence relations, of Eqs. (\ref {17}) and (\ref {19}) the determinants $D(n)$ satisfy the recurrence relations 
\
\begin{eqnarray}
&& D(n)=D(n-3)C(1,1,1)-D(n-4)C(1,1,2)+D(n-5)C(1,1,3)- ...\nonumber\\
&& +(-1)^{n-3}D(2)C(1,1,n-4)+ (-1)^{n-2}D(1)C(1,1,n-3)+ \nonumber\\
&& (-1)^{n-1}C(1,1,n-2),\;n\;\;=4, 5, 6, ...
 \label{20}
\end{eqnarray}
\begin{eqnarray}
&& D(n)=D(n-4)C(1,1,1,1)-D(n-5)C(1,1,1,2)+D(n-6)C(1,1,1,3)\nonumber\\
&&...+(-1)^{n-2}D(2)C(1,1,1,n-5)+ (-1)^{n-1}D(1)C(1,1,1,n-4)+ \nonumber\\
&& (-1)^nC(1,1,1,n-3) \nonumber\\
&& \;n\;\;= 5, 6,7,... \nonumber\\
&& . \;\;\; . \;\;\;.
\label{20}
\end{eqnarray}  
\begin{eqnarray}
&& D(n)=D(1)C(a_1, a_2,...a_{n-1})-C(a_1,a_2,...a_{n-2},2) \nonumber\\
&& a_j=1,\;\;\; j=1, 2,3,...n-1
\label{22}
\end{eqnarray}

Eqs (\ref{8}), (\ref{9}), (\ref{13}) and the recurrence relations  
can be used to express D(n) as a polynomial of the coefficients $C(j)$. We get  for $n=6,\; 7$ and $ 8$

\begin{eqnarray}
&& D(6)=C(1)^6 -5 C(1)^4 C(2) + 6 C(1)^2 C(2)^2-C^3(2) + 4 C(1)^3 C(3) - \nonumber\\
&&  6 C(1) C(2) C(3) + C(3)^2 - 3 C(1)^2 C(4) + 2 C(2) C(4) + \nonumber\\
&&  2 C(1) C(5) - C(6)
\label{22}
\end{eqnarray}

\begin{eqnarray}
&& D(7)=C(1)^7 - 6 C(1)^5 C(2) + 10 C(1)^3 C(2)^2 - 4 C(1) C(2)^3 + 
 5 C(1)^4 C(3)  \nonumber\\
 && - 12 C(1)^2 C(2) C(3) + 3 C(2)^2 C(3) + 
 3 C(1) C(3)^2 - 4 C(1)^3 C(4) + \nonumber\\
 && 6 C(1) C(2) C(4) - 2 C(3) C(4) + 
 3 C(1)^2 C(5) - 2 C(2) C(5) - \nonumber\\
 && 2 C(1) C(6) + C(7) 
 \label{23}
 \end{eqnarray}

 \begin{eqnarray}
&& D(8)=C(1)^8 - 7 C(1)^6 C(2) + 15 C(1)^4 C(2)^2 - 10 C(1)^2 C(2)^3 + C(2)^4 +\nonumber\\
 &&  6 C(1)^5 C(3) - 20 C(1)^3 C(2) C(3) + 12 C(1) C(2)^2 C(3) + 
 6 C(1)^2 C(3)^2 - \nonumber\\
 && 3 C(2) C(3)^2 - 5 C(1)^4 C(4) + 
 12 C(1)^2 C(2) C(4) - 3 C(2)^2 C(4) - \nonumber\\
 && 6 C(1) C(3) C(4) + C(4)^2 + 
 4 C(1)^3 C(5) - 6 C(1) C(2) C(5) + 2 C(3) C(5)  \nonumber\\
 && -3 C(1)^2 C(6) + 
 2 C(2) C(6) + 2 C(1) C(7) -C(8) 
 \label{24}
 \end{eqnarray}


\begin{thebibliography}{99}
\bibitem
{IN1} Inselberg A. On determinants of Toeplitz-Hessenberg matrices arising in power series, J Math. Anal. Appl. 63 (1978), 347
\bibitem 
{ME1} Merca M A note on the determinant of a Toeplitz- 
Hessenberg maqtrix Special Matrices 1 (2013) 10
\bibitem
{GO1} Goy T and Shattuck M Determinants of Hessenberg - Toeplitz matrises with generalized Fibonacci entries Notes on Number Theory and Discrete Mathematics 25 (2019) 83 

\end{thebibliography}
\end{document}